\documentclass[reqno,a4paper]{amsart}
\usepackage{amsthm}
\usepackage[T1]{fontenc}
\usepackage[utf8]{inputenc}
\usepackage[english]{babel}
\usepackage{microtype}

\usepackage{amssymb, amsmath}
\usepackage{enumitem}
\usepackage{tikz-cd}
\usepackage{orcidlink}
\usepackage{doi}
\renewcommand{\doi}[1]{DOI:~\href{https://doi.org/#1}{#1}}
\usepackage{url}
\usepackage{mathtools}
\usepackage{xurl}
\usepackage{nameref,hyperref,cleveref}
\hypersetup{hidelinks}

\theoremstyle{plain}
\newtheorem{theorem}[subsection]{Theorem}
\newtheorem{proposition}[subsection]{Proposition}
\newtheorem{lemma}[subsection]{Lemma}
\newtheorem{corollary}[subsection]{Corollary}

\theoremstyle{definition}
\newtheorem{definition}[subsection]{Definition}

\theoremstyle{remark}
\newtheorem{examples}[subsection]{Examples}

\crefname{examples}{Examples}{Examples}
\Crefname{examples}{Examples}{Examples}

\numberwithin{equation}{section}

\newcommand{\two}{\mathbf{2}}
\newcommand{\GHoops}{\ensuremath{\mathsf{GHoops}}}
\newcommand{\GAlg}{\ensuremath{\mathsf{GAlg}}}
\newcommand{\BHoops}{\ensuremath{\mathsf{BHoops}}}
\newcommand{\BLAlg}{\ensuremath{\mathsf{BLAlg}}}
\newcommand{\one}{\mathbf{1}}
\newcommand{\up}{{\uparrow}}
\newcommand{\down}{{\downarrow}}
\DeclareMathOperator{\id}{id}

\begin{document}

\title[Free G\"odel algebras over G\"odel hoops]{Free G\"odel algebras over G\"odel hoops}

\author[M.~Mancini]{Manuel Mancini~\orcidlink{0000-0003-2142-6193}}

\email{manuel.mancini@uclouvain.be; manuel.mancini@unipa.it}

\address{Institut de Recherche en Mathématique et Physique, Université catholique de Louvain, chemin du cyclotron 2 bte L7.01.02, B--1348 Louvain-la-Neuve, Belgium.}
\address{Dipartimento di Matematica e Informatica, Università degli Studi di Palermo, Via Archirafi 34, 90123 Palermo, Italy.}

\date{}

\begin{abstract}
We give an explicit description of the left adjoint $F$ of the forgetful functor $U\colon\GAlg\to\GHoops$ from G\"odel algebras to G\"odel hoops. If $H$ is a G\"odel hoop with a least element~$m$, then $F(H)$ is the direct product of the G\"odel algebra obtained by adjoining a new bottom element to~$H$ and of $H$ itself, regarded as a G\"odel algebra with bottom element~$m$. For an arbitrary G\"odel hoop~$H$, we glue these products along the principal filters of~$H$, and we describe the operations of $F(H)$, the action of $F$ on homomorphisms and the unit of the adjunction. As an application, we show that the unit of the adjunction $F\dashv U$ is not cartesian, and that the same holds for BL-algebras and basic hoops. Moreover, the semi-abelian categories of G\"odel hoops and of basic hoops are not equivalent to the slice categories of G\"odel algebras and of BL-algebras over their initial object. Thus, in contrast with the cases of Wajsberg hoops and of product hoops, the free--forgetful adjunctions between G\"odel algebras and G\"odel hoops, and between BL-algebras and basic hoops, do not give rise to ideally exact contexts.
\end{abstract}

\subjclass[2020]{03G25; 06D20; 06F05; 08B20; 18A40; 18E13}
\keywords{G\"odel hoop, G\"odel algebra, BL-algebra, basic hoop, free--forgetful adjunction, ideally exact category, semi-abelian category}

\maketitle

\section{Introduction}\label{SecIntro}

Among the many-valued logics studied in mathematical fuzzy logic, a central role is played by those whose truth values are the real numbers of the unit interval $[0,1]$, and in which the conjunction is interpreted by a continuous t-norm. Recall that a \emph{t-norm} is a binary operation $\ast$ on $[0,1]$ which is associative, commutative, monotone in each argument with respect to the usual order, and has $1$ as neutral element. If $\ast$ is continuous, then the implication is interpreted by its residuum $x\to y=\max\{z\in[0,1]\mid z\ast x\leq y\}$. The three fundamental continuous t-norms are the \L{}ukasiewicz t-norm $x\ast y=\max\{x+y-1,0\}$, the G\"odel t-norm $x\ast y=\min\{x,y\}$ and the product t-norm $x\ast y=xy$. Every continuous t-norm is an ordinal sum of isomorphic copies of the \L{}ukasiewicz and product t-norms, defined on pairwise disjoint open subintervals of $[0,1]$, and it coincides with the G\"odel t-norm whenever its two arguments do not lie in the same subinterval~\cite{Hajek}. The corresponding logics are \L{}ukasiewicz infinite-valued logic, G\"odel logic and product logic, whose equivalent algebraic semantics are given by the varieties of MV-algebras~\cite{Chang}, G\"odel algebras~\cite{Hajek} and product algebras~\cite{HGE}, respectively. In~\cite{Hajek}, P.~H\'ajek introduced Basic Logic, a common fragment of these three logics, together with its algebraic semantics given by BL-algebras, and he conjectured that it is the logic of all continuous t-norms and their residua. This conjecture was proved in~\cite{CEGT}, and a more algebraic proof was given in~\cite{AFM}.

Hoops first appeared in the work of B.~Bosbach on \emph{complementary semigroups}~\cite{Bosbach}, and their structure theory was developed by W.~J.~Blok and I.~M.~A.~Ferreirim~\cite{BlokFerreirim}. Basic hoops, Wajsberg hoops, G\"odel hoops and product hoops are the $0$-free subreducts of BL-algebras, MV-algebras, G\"odel algebras and product algebras, respectively~\cite{AFM,GMU}, and the latter are term equivalent to the corresponding classes of bounded hoops, that is, of hoops with a least element named by a constant~$0$; see~\cite{AFM,Hajek}. In each case, the functor which forgets the constant~$0$ has a left adjoint, which freely adjoins a falsum constant. An explicit description of this left adjoint is known for Wajsberg hoops, where it is given by the \emph{MV-closure} of~\cite{ACD}, and for product hoops, where it was obtained in~\cite{GMU} by means of $\ell$-groups. On the other hand, it was observed in~\cite[Section~6]{GMU} that it is not known how to freely obtain a G\"odel algebra from a G\"odel hoop, and that the methods developed there for product hoops cannot be directly applied to this case.

In this article, we give such a construction. Let $U\colon\GAlg\to\GHoops$ be the forgetful functor from G\"odel algebras to G\"odel hoops, and let $F$ be its left adjoint. For a G\"odel hoop~$H$, we denote by $H^+$ the G\"odel algebra obtained by adjoining a new bottom element to~$H$. We show in \Cref{LemExtension} that, if $H$ has a least element~$m$, then $F(H)$ is the direct product of $H^+$ and of $H$ itself, regarded as a G\"odel algebra with bottom element~$m$. Adjoining a new bottom element is not enough, since a homomorphism of G\"odel hoops from $H$ to a G\"odel algebra may send $m$ to~$0$: for homomorphisms into a G\"odel chain, the factor~$H^+$ accounts for those which send $m$ to a non-zero element, and the factor~$H$ for those which send $m$ to~$0$. Since every G\"odel hoop is the directed union of its principal filters, each of which has a least element, $F(H)$ can be obtained in general by gluing together these products. In \Cref{ThmMain}, we make this explicit: we describe $F(H)$ as a quotient of a set of triples of elements of~$H^+$, with operations given by explicit terms in the language of~$H$, and we verify directly that $F(H)$ is a G\"odel algebra with the required universal property. For instance, \Cref{ExExamples} show that $F(H)\cong H^+$ whenever $H$ is a totally ordered G\"odel hoop without least element. The G\"odel algebra $H^+$ also appears in~\cite[Section~3.1]{FG}, in connection with Esakia duality for relative Stone algebras, while G\"odel algebras free over distributive lattices, which satisfy a different universal property, were described via Esakia duality in~\cite[Theorem~3.19]{Carai}; see also~\cite{AGM} for the finite case. As a consequence, we recover the result of P.~K\"ohler which expresses finitely generated free G\"odel algebras in terms of free G\"odel hoops; see \Cref{CorFree} and~\cite[Theorem~5.7]{Davey}.

Our motivation comes from the study of internal actions in ideally exact categories. The variety $\GHoops$ is a semi-abelian category~\cite{JMT,LMS}, while $\GAlg$ is an \emph{ideally exact category} in the sense of G.~Janelidze~\cite{IdeallyExact}. By~\cite{IdeallyExact} and~\cite[Remark~3.3]{MMP}, a monadic functor from an ideally exact category $\mathcal{U}$, with initial object~$0$ and terminal object~$1$, to a semi-abelian category is, up to an equivalence, the pullback functor $\mathcal{U}\to(\mathcal{U}\downarrow0)$ along the unique morphism $0\to1$ if and only if the unit of the corresponding adjunction is cartesian. Monadic adjunctions of this kind, called \emph{ideally exact contexts}, were used in~\cite{MMP} to introduce coherent and ideal actions, and it was shown in~\cite[Theorems~6.30 and~6.39]{MMP} that MV-algebras and Wajsberg hoops, as well as product algebras and product hoops, form ideally exact contexts with a \emph{good theory of actions}, that is, in which every coherent action is ideal, and so is every morphism between such actions. Internal actions and split extensions of hoops were also studied in~\cite{MMPStudia,MMPT}. In~\cite{MP}, the ideally exact contexts with a good theory of actions were characterised among those obtained from a semi-abelian variety by freely adding constants, provided that the unit of the corresponding free--forgetful adjunction is cartesian. As observed in~\cite[Section~4]{MP}, no explicit description was available for the left adjoints of the forgetful functors $\GAlg\to\GHoops$ and $\BLAlg\to\BHoops$, where $\BLAlg$ and $\BHoops$ denote the varieties of BL-algebras and of basic hoops, and it was not known whether the units of the corresponding adjunctions are cartesian.

We show that neither of these units is cartesian, so that the corresponding free--forgetful adjunctions are not ideally exact contexts. For G\"odel hoops, this is proved in \Cref{PropNotCartesian}, and the case of basic hoops follows in \Cref{CorBL}. Moreover, by \Cref{PropNotEquivalent}, the semi-abelian categories $\GHoops$ and $(\GAlg\downarrow\two)$, where $\two$ denotes the two-element Boolean algebra, are not equivalent, and, by \Cref{PropBLNotEquivalent}, neither are $\BHoops$ and $(\BLAlg\downarrow\two)$. This corrects the claim, made in the introduction of~\cite{MMPStudia} but not used elsewhere in that paper, that these categories are equivalent.

\section{Preliminaries}\label{SecPrelim}

\begin{definition}\label{DefGodel}
A \emph{G\"odel hoop} is an algebra $(H,\wedge,\vee,\to,1)$ such that $(H,\wedge,\vee,1)$ is a lattice with top element~$1$ and, for all $x,y,z\in H$,
\begin{align}
x\wedge y\leq z\quad&\Longleftrightarrow\quad x\leq y\to z,\label{EqResiduation}\\
(x\to y)\vee(y\to x)&=1.\label{EqPrelinearity}
\end{align}
A \emph{G\"odel algebra} is an algebra $(A,\wedge,\vee,\to,0,1)$ such that $(A,\wedge,\vee,\to,1)$ is a G\"odel hoop and $0\leq x$ for every $x\in A$.
\end{definition}

Thus, G\"odel hoops are the prelinear Brouwerian algebras, also known as relative Stone algebras~\cite{Davey,FG}, and G\"odel algebras are the prelinear Heyting algebras. G\"odel hoops are usually defined as the basic hoops $(H,\cdot,\to,1)$ satisfying $x\cdot x=x$~\cite{AFM}: the two presentations are term equivalent, since in such hoops $x\cdot y=x\wedge y$ and $x\vee y=((x\to y)\to y)\wedge((y\to x)\to x)$. G\"odel hoops and G\"odel algebras form varieties, which we denote by $\GHoops$ and $\GAlg$. Homomorphisms preserve all the operations and constants of the respective signatures. In particular, a homomorphism of G\"odel hoops between two G\"odel algebras need not preserve their least elements. The underlying lattice of a G\"odel hoop is distributive, and we freely use elementary identities such as $x\wedge(x\to y)=x\wedge y$, $y\leq x\to y$ and $x\to(y\to z)=(x\wedge y)\to z$, as well as the fact that $\to$ is antitone in its first and monotone in its second argument.

The forgetful functor $U\colon\GAlg\to\GHoops$, which forgets the constant~$0$, has a left adjoint~$F$, since it is induced by an inclusion of signatures between varieties. In a G\"odel algebra we write $\neg x=x\to0$. We denote by $\two=\{0<1\}$ the two-element Boolean algebra, which is the initial G\"odel algebra, by $\one=\{1\}$ the trivial G\"odel hoop, which is the zero object of $\GHoops$, and by $C_n$ the $n$-element chain, regarded as a G\"odel algebra or as a G\"odel hoop.

\begin{lemma}\label{LemStone}
Let $A$ be a G\"odel algebra.
\begin{enumerate}
\item For every $e\in A$, the interval ${\down}e=\{x\in A\mid x\leq e\}$, with the restricted lattice operations, the bottom element~$0$, the top element~$e$ and the implication $x\to_e y=e\wedge(x\to y)$, is a G\"odel algebra, and the map $p_e\colon A\to{\down}e\colon x\mapsto x\wedge e$ is a surjective homomorphism of G\"odel algebras.
\item For every $a\in A$, the elements $c=\neg a$ and $d=\neg\neg a$ satisfy $c\wedge d=0$ and $c\vee d=1$, and the map $A\to{\down}d\times{\down}c\colon x\mapsto(x\wedge d,x\wedge c)$ is an isomorphism of G\"odel algebras, whose inverse is $(u,v)\mapsto u\vee v$.
\end{enumerate}
\end{lemma}

\begin{proof}
(1) For all $x,y\in A$ and $z\in{\down}e$, we have
\[
z\leq p_e(x\to y)\iff z\wedge x\leq y\iff z\wedge x\wedge e\leq y\wedge e\iff z\leq p_e(x)\to_e p_e(y),
\]
since $z\wedge x\leq e$. If $x,y\in{\down}e$, then $p_e(x)=x$ and $p_e(y)=y$, and these equivalences show that ${\down}e$ satisfies~\eqref{EqResiduation} with respect to~$\to_e$. For arbitrary $x,y\in A$, they show that $p_e(x\to y)=p_e(x)\to_e p_e(y)$. Clearly, $p_e$ preserves meets and constants, and it preserves joins by distributivity. Since $p_e$ is surjective, \eqref{EqPrelinearity} holds in~${\down}e$.

(2) We have $c\wedge d=c\wedge(c\to0)=c\wedge0=0$ and $a\to c=(a\wedge a)\to0=c$. Since $a\wedge c=0$, for every $z\in A$ we have $z\wedge c\leq a$ if and only if $z\wedge c=0$, so that $c\to a=c\to0=d$. Hence, \eqref{EqPrelinearity} applied to $a$ and $c$ gives $c\vee d=1$. The map in~(2) is a homomorphism by~(1), and it is bijective, with the stated inverse, since $c$ and $d$ are complementary elements of a distributive lattice.
\end{proof}

For a G\"odel hoop~$H$, we denote by $H^+=H\sqcup\{\bot\}$ the set obtained by adding to~$H$ a new element~$\bot$, declared to be strictly below every element of~$H$. We keep the operations of~$H$, and for $u\in H^+$ and $x\in H$ we put
\[
\bot\wedge u=u\wedge\bot=\bot,\quad \bot\vee u=u\vee\bot=u,\quad \bot\to u=1,\quad x\to\bot=\bot.
\]

\begin{lemma}\label{LemPlus}
Let $H$ be a G\"odel hoop.
\begin{enumerate}
\item $H^+$ is a G\"odel algebra with bottom element~$\bot$.
\item Every homomorphism of G\"odel hoops $f\colon H\to K$ extends to a homomorphism of G\"odel algebras $f^+\colon H^+\to K^+$ with $f^+(\bot)=\bot$.
\item Let $A$ be a G\"odel algebra and let $g\colon H\to U(A)$ be a homomorphism of G\"odel hoops. Then the extension $g^\sharp\colon H^+\to A$ of $g$ with $g^\sharp(\bot)=0$ is a homomorphism of G\"odel algebras if and only if $\neg g(x)=0$ for every $x\in H$.
\end{enumerate}
\end{lemma}

\begin{proof}
(1) Condition~\eqref{EqPrelinearity} holds whenever one of its arguments is~$\bot$, since $\bot\to u=1$ for every $u\in H^+$. As for~\eqref{EqResiduation}, both sides hold if $y=\bot$, both sides are equivalent to $x=\bot$ if $y\in H$ and $z=\bot$, and in the remaining cases \eqref{EqResiduation} holds in $H$ or trivially.

(2) This follows from the definitions, since $f$ maps $H$ into~$K$.

(3) The map $g^\sharp$ preserves the constants and all the operations, except possibly the implications $x\to\bot=\bot$ with $x\in H$, which are preserved if and only if $g(x)\to0=0$, that is, $\neg g(x)=0$.
\end{proof}

Finally, let $H$ be a G\"odel hoop and let $b\in H$. The principal filter ${\up}b$ of the elements $x\in H$ such that $b\leq x$ is closed under the operations of~$H$, since $y\leq x\to y$. We denote by $H_b$ the G\"odel algebra obtained by regarding ${\up}b$ as a G\"odel hoop with bottom element~$b$, and we put $H_b^+=({\up}b)^+$. Notice that the operations $\wedge$, $\vee$ and $\to$ of the G\"odel algebra $H_b^+\times H_b$ are the restrictions of those of $H^+\times H$, that its bottom element is $(\bot,b)$, and that pseudocomplements in $H_b$ are given by $\neg_b x=x\to b$.

\section{The construction of the free functor}\label{SecFree}

Let $H$ be a G\"odel hoop. For every $b\in H$, the principal filter ${\up}b$ is a G\"odel hoop with least element~$b$. We first show in \Cref{LemExtension} that $H_b^+\times H_b$ is the free G\"odel algebra on~${\up}b$, and we describe in \Cref{LemTransition} the homomorphisms between these free algebras induced by the inclusions ${\up}b\subseteq{\up}c$, for $c\leq b$. Since $H$ is the directed union of its principal filters, $F(H)$ is then obtained in \Cref{ThmMain} by gluing the algebras $H_b^+\times H_b$ along these homomorphisms.

\begin{lemma}\label{LemExtension}
Let $H$ be a G\"odel hoop, let $b\in H$, let $A$ be a G\"odel algebra and let $g\colon{\up}b\to U(A)$ be a homomorphism of G\"odel hoops. Then the map $\bar g\colon H_b^+\times H_b\to A$ defined by
\begin{equation}\label{EqExtension}
\bar g(u,v)=\bigl(g^\sharp(u)\wedge\neg\neg g(b)\bigr)\vee\bigl(g(v)\wedge\neg g(b)\bigr)
\end{equation}
is the unique homomorphism of G\"odel algebras such that $\bar g(x,x)=g(x)$ for every $x\in{\up}b$.
\end{lemma}

\begin{proof}
Put $c=\neg g(b)$ and $d=\neg\neg g(b)$. By \Cref{LemStone}(1), $p_c\circ g\colon{\up}b\to{\down}c$ is a homomorphism of G\"odel hoops which sends $b$ to $g(b)\wedge c=0$, hence a homomorphism of G\"odel algebras $H_b\to{\down}c$. For every $x\in{\up}b$ we have $g(b)\leq g(x)\wedge d$, so that the pseudocomplement $d\wedge\neg(g(x)\wedge d)$ of $p_d(g(x))$ in ${\down}d$ is below $d\wedge c=0$. Thus, by \Cref{LemPlus}(3), the extension $(p_d\circ g)^\sharp\colon H_b^+\to{\down}d$ is a homomorphism of G\"odel algebras. Composing $(p_d\circ g)^\sharp\times(p_c\circ g)$ with the inverse of the isomorphism of \Cref{LemStone}(2), we obtain that $\bar g$ is a homomorphism of G\"odel algebras. Moreover, $\bar g(x,x)=g(x)\wedge(d\vee c)=g(x)$ for every $x\in{\up}b$.

For the uniqueness, notice that the bottom element of $H_b^+\times H_b$ is $(\bot,b)$, so that $\neg(b,b)=(\bot,1)$ and $\neg\neg(b,b)=(1,b)$. Hence, for all $u,v\in{\up}b$,
\begin{equation}\label{EqGenerators}
(u,v)=\bigl((u,u)\wedge\neg\neg(b,b)\bigr)\vee\bigl((v,v)\wedge\neg(b,b)\bigr),\qquad(\bot,v)=(v,v)\wedge\neg(b,b).
\end{equation}
Therefore, $H_b^+\times H_b$ is generated, as a G\"odel algebra, by its diagonal $\{(x,x)\mid x\in{\up}b\}$, and the uniqueness follows.
\end{proof}

In other words, $H_b^+\times H_b$, together with the map $x\mapsto(x,x)$, is the free G\"odel algebra on the G\"odel hoop ${\up}b$. The definition of~$\bar g$ reflects the decomposition $A\cong{\down}\neg\neg g(b)\times{\down}\neg g(b)$ of \Cref{LemStone}(2): in the first factor, the image of every $x\in{\up}b$ has pseudocomplement~$0$, so that $\bot$ can be sent to the bottom element, while in the second factor $g(b)$ becomes the bottom element, so that $H_b$ is mapped there as a G\"odel algebra.

Let $c\leq b$ in~$H$, so that $\neg_c b=b\to c$ and $\neg_c\neg_c b=(b\to c)\to c$ are the pseudocomplement and the double pseudocomplement of $b$ in~$H_c$. For $(u,v)\in H_b^+\times H_b$, we put
\begin{equation}\label{EqRho}
\rho_{b,c}(u,v)=(u\wedge\neg_c\neg_c b)\vee(v\wedge\neg_c b),
\end{equation}
where the operations are computed in~$H^+$. Since $c\leq v\wedge\neg_c b$, we have $\rho_{b,c}(u,v)\in{\up}c$. For instance, $\rho_{b,c}(\bot,v)=v\wedge\neg_c b$. With respect to the decomposition $H_c\cong{\down}\neg_c\neg_c b\times{\down}\neg_c b$ of \Cref{LemStone}(2), the first component of $\rho_{b,c}(u,v)$ only depends on~$u$, and the second one only on~$v$.

\begin{lemma}\label{LemTransition}
Let $H$ be a G\"odel hoop and let $b,c\in H$ with $c\leq b$. Then the map $\gamma_{b,c}\colon H_b^+\times H_b\to H_c^+\times H_c$ defined by
\[
\gamma_{b,c}(u,v)=\bigl(u,\rho_{b,c}(u,v)\bigr)
\]
is the unique homomorphism of G\"odel algebras such that $\gamma_{b,c}(x,x)=(x,x)$ for every $x\in{\up}b$. Moreover, $\gamma_{b,b}$ is the identity and $\gamma_{c,d}\circ\gamma_{b,c}=\gamma_{b,d}$ whenever $d\leq c\leq b$.
\end{lemma}

\begin{proof}
We apply \Cref{LemExtension} to the homomorphism of G\"odel hoops $g\colon{\up}b\to U(H_c^+\times H_c)\colon x\mapsto(x,x)$. Write $p=\neg_c b$ and $q=\neg_c\neg_c b$. Since the bottom element of $H_c^+\times H_c$ is $(\bot,c)$, we have $g^\sharp(\bot)=(\bot,c)$ and
\[
\neg g(b)=(b\to\bot,b\to c)=(\bot,p),\qquad \neg\neg g(b)=(\bot\to\bot,p\to c)=(1,q).
\]
Hence, by~\eqref{EqExtension}, we get
\[
\bar g(u,v)=\bigl(u,(u\wedge q)\vee(v\wedge p)\bigr)\quad\text{if }u\in{\up}b,\qquad \bar g(\bot,v)=\bigl(\bot,(c\wedge q)\vee(v\wedge p)\bigr).
\]
Since $c\wedge q=c\leq v\wedge p$, in both cases we have $\bar g(u,v)=\gamma_{b,c}(u,v)$. The last assertion follows from the uniqueness in \Cref{LemExtension}, since $\gamma_{b,b}$, $\gamma_{c,d}\circ\gamma_{b,c}$ and $\gamma_{b,d}$ send $(x,x)$ to $(x,x)$ for every $x\in{\up}b$.
\end{proof}

We are now ready to describe $F(H)$. Let
\[
T(H)=\bigl\{(b,u,v)\in H\times H^+\times H\mid (u,v)\in H_b^+\times H_b\bigr\},
\]
that is, $(b,u,v)\in T(H)$ if and only if $b\leq v$ and either $u=\bot$ or $b\leq u$. We write $(b,u,v)\sim(b',u',v')$ if there exists $c\leq b\wedge b'$ such that $\gamma_{b,c}(u,v)=\gamma_{b',c}(u',v')$, that is, such that $u=u'$ and $\rho_{b,c}(u,v)=\rho_{b',c}(u',v')$.

Thus, a triple $(b,u,v)$ represents the element $(u,v)$ of the free G\"odel algebra $H_b^+\times H_b$ on ${\up}b$, and two triples are equivalent when they have the same image at some common stage $c\leq b\wedge b'$. Since $H$ is the directed union of its principal filters and $F$, being a left adjoint, preserves directed colimits, $F(H)$ is the directed colimit of the G\"odel algebras $H_b^+\times H_b$ along the homomorphisms~$\gamma_{b,c}$. In the following theorem, we verify directly that the quotient $T(H)/{\sim}$ is a G\"odel algebra with the required universal property, and we make its structure explicit.

\begin{theorem}\label{ThmMain}
Let $H$ be a G\"odel hoop.
\begin{enumerate}
\item The relation $\sim$ is an equivalence relation on $T(H)$. If $[b,u,v]$ denotes the class of $(b,u,v)$, then $F(H)=T(H)/{\sim}$ is a G\"odel algebra with constants $0=[b,\bot,b]$ and $1=[b,1,1]$, which do not depend on $b\in H$, and with operations
\begin{equation}\label{EqOperations}
[b,u,v]\star[b',u',v']=\bigl[c,\ u\star u',\ \rho_{b,c}(u,v)\star\rho_{b',c}(u',v')\bigr]
\end{equation}
for $\star\in\{\wedge,\vee,\to\}$, where $c=b\wedge b'$ and the operations in the second and in the third entry are computed in~$H^+$ and in~$H$, respectively.
\item The map $\eta_H\colon H\to UF(H)\colon h\mapsto[h,h,h]$ is a homomorphism of G\"odel hoops. For every G\"odel algebra $A$ and every homomorphism of G\"odel hoops $g\colon H\to U(A)$, the map $\hat g\colon F(H)\to A$ defined by
\begin{equation}\label{EqHat}
\hat g[b,u,v]=\bigl(g^\sharp(u)\wedge\neg\neg g(b)\bigr)\vee\bigl(g(v)\wedge\neg g(b)\bigr)
\end{equation}
is the unique homomorphism of G\"odel algebras such that $U(\hat g)\circ\eta_H=g$. Hence, $F$ is the left adjoint of~$U$, with unit~$\eta$.
\item For every homomorphism of G\"odel hoops $f\colon H\to K$, the homomorphism $F(f)\colon F(H)\to F(K)$ is given by
\[
F(f)[b,u,v]=\bigl[f(b),f^+(u),f(v)\bigr].
\]
\end{enumerate}
\end{theorem}

\begin{proof}
For $(b,u,v)\in T(H)$, we write $X=(u,v)$ and $[b,X]=[b,u,v]$, so that~\eqref{EqOperations} reads
\[
[b,X]\star[b',X']=\bigl[c,\gamma_{b,c}(X)\star\gamma_{b',c}(X')\bigr],
\]
where the operation on the right-hand side is computed in the G\"odel algebra $H_c^+\times H_c$.

(1) The relation $\sim$ is clearly reflexive and symmetric. If $(b,X)\sim(b',X')$ and $(b',X')\sim(b'',X'')$, with witnesses $c$ and $c'$ respectively, then $d=c\wedge c'$ satisfies
\begin{align*}
\gamma_{b,d}(X)&=\gamma_{c,d}\gamma_{b,c}(X)=\gamma_{c,d}\gamma_{b',c}(X')=\gamma_{b',d}(X')\\
&=\gamma_{c',d}\gamma_{b',c'}(X')=\gamma_{c',d}\gamma_{b'',c'}(X'')=\gamma_{b'',d}(X'')
\end{align*}
by \Cref{LemTransition}, so that $\sim$ is transitive. Moreover, $(b,X)\sim(c,\gamma_{b,c}(X))$ for every $c\leq b$, that is,
\begin{equation}\label{EqStage}
[b,X]=\bigl[c,\gamma_{b,c}(X)\bigr]\qquad\text{whenever }c\leq b.
\end{equation}

Next, we show that the operations~\eqref{EqOperations} are well defined. Let $d\leq c=b\wedge b'$. Since $\gamma_{c,d}$ is a homomorphism, \eqref{EqStage} and \Cref{LemTransition} give
\[
\bigl[c,\gamma_{b,c}(X)\star\gamma_{b',c}(X')\bigr]=\bigl[d,\gamma_{c,d}\bigl(\gamma_{b,c}(X)\star\gamma_{b',c}(X')\bigr)\bigr]=\bigl[d,\gamma_{b,d}(X)\star\gamma_{b',d}(X')\bigr],
\]
so that the right-hand side of~\eqref{EqOperations} may be computed at any $d\leq b\wedge b'$. Now, if $(b,X)\sim(b_1,X_1)$ with witness $e$, and $d=e\wedge b'$, then $\gamma_{b,d}(X)=\gamma_{e,d}\gamma_{b,e}(X)=\gamma_{e,d}\gamma_{b_1,e}(X_1)=\gamma_{b_1,d}(X_1)$, so that $[b,X]\star[b',X']=[b_1,X_1]\star[b',X']$. The same argument applies to the second variable. Similarly, the constants do not depend on~$b$, since each $\gamma_{b,c}$ sends the bottom element $(\bot,b)$ and the top element $(1,1)$ of $H_b^+\times H_b$ to those of $H_c^+\times H_c$.

It remains to show that $F(H)$ is a G\"odel algebra. For $c\in H$, let $S_c$ be the image of the map
\[
\lambda_c\colon H_c^+\times H_c\to F(H)\colon X\mapsto[c,X].
\]
By~\eqref{EqOperations} and the definition of $0$ and~$1$, $\lambda_c$ preserves the operations and the constants. Hence, $S_c$ is closed under the operations of $F(H)$, and it is a homomorphic image of the G\"odel algebra $H_c^+\times H_c$. Since G\"odel algebras form a variety, $S_c$ is a G\"odel algebra. Moreover, $F(H)$ is the union of the subsets $S_b$, with $b\in H$, and $S_b\subseteq S_c$ whenever $c\leq b$, by~\eqref{EqStage}. Thus, any finitely many elements of $F(H)$ belong to a common $S_c$, namely for $c$ the meet of the first components of chosen representatives. Since each condition in \Cref{DefGodel} involves finitely many elements, $F(H)$ is a G\"odel algebra.

(2) Let $x,y\in H$, and put $c=x\wedge y$ and $z=x\star y$, so that $c\leq z$ (if $\star$ is~$\to$, since $c\leq y\leq x\to y$). By \Cref{LemTransition} and~\eqref{EqStage}, we have
\[
\eta_H(x)\star\eta_H(y)=\bigl[c,(x,x)\star(y,y)\bigr]=[c,z,z]=[z,z,z]=\eta_H(z).
\]
Moreover, $\eta_H(1)=[1,1,1]=1$, so that $\eta_H$ is a homomorphism of G\"odel hoops.

Let $g\colon H\to U(A)$ be a homomorphism of G\"odel hoops and, for $b\in H$, let $\bar g_b\colon H_b^+\times H_b\to A$ be the homomorphism of \Cref{LemExtension} associated with the restriction of $g$ to~${\up}b$. If $c\leq b$, then $\bar g_c\circ\gamma_{b,c}$ and $\bar g_b$ are homomorphisms of G\"odel algebras sending $(x,x)$ to $g(x)$ for every $x\in{\up}b$, so that $\bar g_c\circ\gamma_{b,c}=\bar g_b$ by \Cref{LemExtension}. It follows that $\hat g[b,X]=\bar g_b(X)$ is well defined, and that $\hat g\circ\lambda_c=\bar g_c$ for every $c\in H$. Since any two elements of $F(H)$ belong to a common $S_c$, and $\lambda_c$ and $\bar g_c$ are homomorphisms, so is~$\hat g$. Moreover, $\hat g(\eta_H(h))=\bar g_h(h,h)=g(h)$ for every $h\in H$.

For the uniqueness, notice that $\lambda_b(x,x)=[b,x,x]=[x,x,x]=\eta_H(x)$ for every $x\in{\up}b$, by~\eqref{EqStage}. Applying the homomorphism $\lambda_b$ to~\eqref{EqGenerators}, we get
\begin{equation}\label{EqEta}
[b,u,v]=\bigl(\eta_H(u)\wedge\neg\neg\eta_H(b)\bigr)\vee\bigl(\eta_H(v)\wedge\neg\eta_H(b)\bigr)
\end{equation}
for every $(b,u,v)\in T(H)$, where $\eta_H(\bot)$ stands for~$0$. Hence, $F(H)$ is generated by $\eta_H(H)$ as a G\"odel algebra, and $\hat g$ is unique.

(3) By~(2), $F(f)$ is the unique homomorphism of G\"odel algebras such that $F(f)\circ\eta_H=\eta_K\circ f$. Applying~\eqref{EqHat} to $g=\eta_K\circ f$, for which $g^\sharp(u)=\eta_K(f^+(u))$ with the convention $\eta_K(\bot)=0$, and then~\eqref{EqEta} in~$F(K)$, we get $F(f)[b,u,v]=[f(b),f^+(u),f(v)]$.
\end{proof}

The homomorphism~$\hat g$ has a particularly simple form when $A$ is a G\"odel chain, since then $\neg g(b)$ is $1$ if $g(b)=0$, and $0$ otherwise. Hence, $\hat g[b,u,v]=g^\sharp(u)$ if $g(b)\neq0$, while $\hat g[b,u,v]=g(v)$ if $g(b)=0$: the third entry~$v$ of a triple $(b,u,v)$ only matters for the homomorphisms which send $b$ to~$0$.

The map $\pi_H\colon F(H)\to H^+\colon[b,u,v]\mapsto u$ is a surjective homomorphism of G\"odel algebras by~\eqref{EqOperations}, and $\pi_H\circ\eta_H$ is the inclusion of $H$ into~$H^+$. In particular, the unit~$\eta_H$ is injective.

\begin{examples}\label{ExExamples}{\ }
\begin{enumerate}
\item Suppose that $H$ has a least element~$m$, and denote by $H_m$ the G\"odel algebra obtained by regarding $H$ as a G\"odel hoop with bottom element~$m$. Then the homomorphism $\lambda_m\colon H^+\times H_m\to F(H)\colon(u,v)\mapsto[m,u,v]$ is an isomorphism. Indeed, it is surjective by~\eqref{EqStage}, and it is injective since $m$ is the only element below~$m$ and $\gamma_{m,m}$ is the identity. Under this isomorphism, $\eta_H(x)$ corresponds to $(x,x)$. For instance, $F(\one)\cong\one^+\times\one\cong\two$ and $F(C_2)\cong C_3\times C_2$. In particular, if $A$ is a G\"odel algebra, then $FU(A)\cong A^+\times A$, and by~\eqref{EqHat} the counit $FU(A)\to A$ of the adjunction is the second projection.
\item Suppose that $H$ is a totally ordered G\"odel hoop without least element. In such a hoop, we have $x\to y=y$ whenever $y<x$. Let $(b,u,v),(b',u,v')\in T(H)$, and choose $c<b\wedge b'$. Then $\neg_c b=c$ and $\neg_c\neg_c b=1$, so that $\rho_{b,c}(u,v)=u\vee(v\wedge c)=u\vee c$, and similarly $\rho_{b',c}(u,v')=u\vee c$. Hence, $[b,u,v]=[b',u,v']$, and the homomorphism~$\pi_H$ is an isomorphism. In other words, $F(H)$ is obtained by adjoining a new bottom element to~$H$. For instance, $F((0,1])$ is the standard G\"odel algebra $[0,1]$.
\item Let $C_2=\{m<1\}$ and let $H=C_2\times(0,1]$, which is not totally ordered and has no least element. We write $x=(x_1,x_2)$ for $x\in H$. By~\eqref{EqStage}, every element of $F(H)$ is of the form $[b,u,v]$ with $b=(m,t)$. If $c=(m,s)$ with $s<t$, then $\neg_c b=(1,s)$ and $\neg_c\neg_c b=(m,1)$, so that $\rho_{b,c}(u,v)$ is $(v_1,u_2)$ if $u\neq\bot$, and $(v_1,s)$ otherwise. Hence, $[b,u,v]$ only depends on $u$ and~$v_1$, and \eqref{EqHat}, applied to $x\mapsto(x,x_1)$, gives an isomorphism $F(H)\to H^+\times\two$ which sends $[b,u,v]$, with $b=(m,t)$, to $(u,v_1)$, where $C_2$ is identified with~$U(\two)$. In particular, $\pi_H$ is not an isomorphism.
\end{enumerate}
\end{examples}

\begin{corollary}\label{CorFree}
Let $P$ be the free G\"odel hoop on a finite set~$X$, and let $m=\bigwedge X$, where $m=1$ if $X=\emptyset$. Then $m$ is the least element of~$P$, and the free G\"odel algebra on~$X$ is isomorphic to $P^+\times P_m$. In particular, if $P$ has $k$ elements, then the free G\"odel algebra on~$X$ has $k(k+1)$ elements.
\end{corollary}

\begin{proof}
The principal filter ${\up}m$ is closed under the operations of~$P$ and contains~$X$, so that ${\up}m=P$. Since $F$ is a left adjoint, $F(P)$ is the free G\"odel algebra on~$X$, and the result follows from \Cref{ExExamples}(1).
\end{proof}

This is the description of finitely generated free G\"odel algebras given in~\cite[Theorem~5.7]{Davey}, where it is attributed to P.~K\"ohler. For instance, the free G\"odel hoops on one and two generators have $2$ and $18$ elements, so that the corresponding free G\"odel algebras have $6$ and $342$ elements; cf.~\cite{Horn}.

\section{A connection with ideally exact categories}\label{SecIdeallyExact}

Recall that a natural transformation is \emph{cartesian} if all its naturality squares are pullbacks.

\begin{proposition}\label{PropNotCartesian}
The unit $\eta\colon 1_{\GHoops}\Rightarrow UF$ of the adjunction $F\dashv U$ is not cartesian.
\end{proposition}

\begin{proof}
Let $H=C_2=\{m<1\}$, let $\tau\colon H\to\one$ be the unique homomorphism, and consider the naturality square
\[
\begin{tikzcd}
H \arrow[r, "\eta_H"] \arrow[d, "\tau"'] & UF(H) \arrow[d, "UF(\tau)"] \\
\one \arrow[r, "\eta_{\one}"'] & UF(\one).
\end{tikzcd}
\]
By \Cref{ThmMain}(3), we have $F(\tau)[b,u,v]=[1,\tau^+(u),1]$, which is $0$ if $u=\bot$ and $1$ otherwise. Since $\eta_{\one}$ selects the top element of $UF(\one)$ and pullbacks in $\GHoops$ are computed on underlying sets, the pullback of $UF(\tau)$ along $\eta_{\one}$ is the set of the elements $[b,u,v]$ of $F(H)$ such that $u\neq\bot$. By \Cref{ExExamples}(1), this set has four elements, namely the classes $[m,u,v]$ with $u,v\in\{m,1\}$, whereas $H$ has two elements. Hence, the square is not a pullback.
\end{proof}

Let $U_{\mathrm{BL}}\colon\BLAlg\to\BHoops$ be the functor which forgets the constant~$0$. As for G\"odel algebras, $U_{\mathrm{BL}}$ has a left adjoint $F_{\mathrm{BL}}$, and we denote by $\eta^{\mathrm{BL}}$ the unit of the adjunction $F_{\mathrm{BL}}\dashv U_{\mathrm{BL}}$. Recall that G\"odel algebras are the BL-algebras satisfying $x\cdot x=x$, in which case $x\cdot y=x\wedge y$, and that G\"odel hoops are the basic hoops satisfying the same identity~\cite{AFM,Hajek}.

\begin{corollary}\label{CorBL}
The unit $\eta^{\mathrm{BL}}\colon 1_{\BHoops}\Rightarrow U_{\mathrm{BL}}F_{\mathrm{BL}}$ of the adjunction $F_{\mathrm{BL}}\dashv U_{\mathrm{BL}}$ is not cartesian.
\end{corollary}

\begin{proof}
Let $H=C_2=\{m<1\}$ and $\tau\colon H\to\one$ be as in the proof of \Cref{PropNotCartesian}, and regard $H$ as a basic hoop and $F(H)$ as a BL-algebra. Let $q\colon F_{\mathrm{BL}}(H)\to F(H)$ be the homomorphism of BL-algebras such that $q\circ\eta^{\mathrm{BL}}_H=\eta_H$. By~\eqref{EqEta}, $F(H)$ is generated by $\eta_H(H)$ as a G\"odel algebra, hence also as a BL-algebra, since $x\wedge y=x\cdot y$ and $x\vee y=((x\to y)\to y)\wedge((y\to x)\to x)$ in~$F(H)$. Thus, $q$ is surjective. Both $F_{\mathrm{BL}}(\one)$ and $F(\one)$ are the initial BL-algebra~$\two$, and $F(\tau)\circ q=F_{\mathrm{BL}}(\tau)$ by the universal property of~$\eta^{\mathrm{BL}}_H$, since both homomorphisms send $\eta^{\mathrm{BL}}_H(x)$ to~$1$ for every $x\in H$. Let $K$ and $K_{\mathrm{BL}}$ be the preimages of~$1$ under $F(\tau)$ and $F_{\mathrm{BL}}(\tau)$, respectively. If $y\in F_{\mathrm{BL}}(H)$ is such that $q(y)\in K$, then $F_{\mathrm{BL}}(\tau)(y)=F(\tau)(q(y))=1$. Hence, $q$ maps $K_{\mathrm{BL}}$ onto~$K$, which has four elements by the proof of \Cref{PropNotCartesian}. Since $K_{\mathrm{BL}}$ is the pullback of $U_{\mathrm{BL}}F_{\mathrm{BL}}(\tau)$ along $\eta^{\mathrm{BL}}_{\one}$, and $H$ has two elements, the naturality square of $\eta^{\mathrm{BL}}$ at~$\tau$ is not a pullback.
\end{proof}

Recall that a variety is a semi-abelian category if and only if it is pointed and protomodular~\cite{BourJan,JMT}, and that a category is \emph{ideally exact} if it is Barr-exact and protomodular, it has finite coproducts, and the unique morphism from its initial object to its terminal object is a regular epimorphism~\cite{IdeallyExact}. The variety $\GAlg$ is protomodular, since the terms witnessing the protomodularity of $\GHoops$~\cite{BourJan,LMS} also witness that of~$\GAlg$. Since the unique morphism from the initial object $\two$ of $\GAlg$ to the trivial G\"odel algebra is surjective, $\GAlg$ is an ideally exact category~\cite{IdeallyExact}, and the slice category $(\GAlg\downarrow\two)$ is semi-abelian; see~\cite[Remark~3.3]{MMP}.

\begin{proposition}\label{PropNotEquivalent}
The semi-abelian categories $\GHoops$ and $(\GAlg\downarrow\two)$ are not equivalent.
\end{proposition}

\begin{proof}
Recall that an object of a semi-abelian category is \emph{simple} if it is not a zero object and every regular epimorphism defined on it is either an isomorphism or a morphism to a zero object. We show that $\GHoops$ has, up to isomorphism, exactly one simple object, whereas $(\GAlg\downarrow\two)$ has at least two.

The quotients of a G\"odel hoop $H$ correspond to its filters, that is, to the upward closed subsets of $H$ containing~$1$ and closed under~$\wedge$~\cite{BlokFerreirim}. Since every principal filter ${\up}x$ is a filter, if $H$ is simple, then ${\up}x=H$ for every $x\in H$ with $x\neq1$. Hence, every $x\neq1$ in~$H$ is the least element of~$H$, and $H\cong C_2$. Conversely, $C_2$ is simple.

The objects of $(\GAlg\downarrow\two)$ are the pairs $(A,p)$, where $p\colon A\to \two$ is a homomorphism, and the zero object is $(\two,\id_{\two})$. The quotients of $(A,p)$ are the surjective homomorphisms $f\colon A\to B$ through which $p$ factors. Since the congruences of a G\"odel algebra correspond to its filters~\cite{Hajek}, the quotients of $(A,p)$ correspond to the filters of $A$ contained in $P=p^{-1}(1)$, the filter $P$ itself giving the zero object. If $P=\{1\}$, then $(A,p)$ is a zero object, since $p(a)=0$ implies $\neg a\in P$, and hence $a=a\wedge\neg a=0$. If $P$ has more than two elements, then it contains some $y\neq1$ which is not its least element, and the principal filter ${\up}y$ lies strictly between $\{1\}$ and~$P$. Hence, $(A,p)$ is simple if and only if $P$ has exactly two elements, and $(C_3,p)$, where $C_3=\{0<x<1\}$ and $p^{-1}(1)=\{x,1\}$, and $(\two\times\two,\pi_1)$, where $\pi_1$ is the first projection, are non-isomorphic simple objects of $(\GAlg\downarrow\two)$.
\end{proof}

The same arguments show that $\BLAlg$ is an ideally exact category with initial object~$\two$, so that $(\BLAlg\downarrow\two)$ is semi-abelian.

\begin{proposition}\label{PropBLNotEquivalent}
The semi-abelian categories $\BHoops$ and $(\BLAlg\downarrow\two)$ are not equivalent.
\end{proposition}

\begin{proof}
We compare the objects with exactly two subobjects, a property which is invariant under equivalences: up to isomorphism, $(\BLAlg\downarrow\two)$ has at least two of them, whereas $\BHoops$ has only one. The subobjects of a pair $(A,p)$ of $(\BLAlg\downarrow\two)$ correspond to the subalgebras of~$A$. Since the only subalgebras of $C_3$ and of $\two\times\two$ are $\{0,1\}$ and the whole algebra, the pairs $(C_3,p)$ and $(\two\times\two,\pi_1)$ of the proof of \Cref{PropNotEquivalent} are two non-isomorphic examples. It remains to show that every basic hoop with exactly two subhoops is isomorphic to~$C_2$.

Let $H$ be such a hoop, so that $H$ is non-trivial and its only subhoops are $\{1\}$ and~$H$. Every filter of~$H$, that is, every upward closed subset containing~$1$ and closed under the product, is a subhoop, since $y\leq x\to y$. Hence, the only filters of~$H$ are $\{1\}$ and~$H$, so that, for all $d,w\in H$ with $d\neq1$, there exists~$n$ such that $d^n\leq w$. If $H$ has an idempotent element $e\neq1$, then $\{e,1\}$ is a subhoop, and $H=\{e,1\}\cong C_2$. Suppose now that $1$ is the only idempotent element of~$H$. Let $a,b\in H$ with $a\neq1$, and put $w=a\cdot b$, $c=a\to w$ and $d=c\to b$. Then $b\leq c$, and the identity $x\cdot(x\to y)=y\cdot(y\to x)$ gives $a\cdot c=w\cdot(w\to a)=w$ and $c\cdot d=b\cdot(b\to c)=b$. Hence, $w\cdot d=a\cdot c\cdot d=w$, so that $w\cdot d^n=w$ for every~$n$. If $d\neq1$, then $d^n\leq w$ for some~$n$, and $w=w\cdot d^n\leq w\cdot w\leq w$ shows that $w$ is idempotent, which is impossible, since $w\leq a<1$. Thus, $d=1$, that is, $a\to a\cdot b=b$. This identity also holds for $a=1$. Finally, let $x\in H$ with $x\neq1$. Since $x^{2j}\to x^{2k}$ is equal to $x^{2(k-j)}$ if $j\leq k$, and to~$1$ otherwise, the set $\{x^{2k}\mid k\geq0\}$ is a subhoop of~$H$, and it is non-trivial, since $x^2\leq x<1$. Therefore, it contains~$x$, that is, $x=x^{2k}$ for some $k\geq1$. But then $x\leq x^2\leq x$, so that $x$ is idempotent, a contradiction.
\end{proof}

The forgetful functor $U$ is monadic by Beck's theorem~\cite[Section~3.3]{BarrWells}, since it reflects isomorphisms and preserves reflexive coequalisers, which in both varieties are computed on underlying sets. Therefore, by~\cite[Theorem~2.6]{IdeallyExact}, \Cref{PropNotEquivalent} gives another proof of \Cref{PropNotCartesian}, and shows more generally that no monadic functor $\GAlg\to\GHoops$ has a left adjoint with cartesian unit. In the terminology of~\cite{MMP}, the adjunction $F\dashv U$ is not an ideally exact context, in contrast with the cases of Wajsberg hoops and product hoops; see~\cite[Section~6.14]{MMP}. By the same argument, $U_{\mathrm{BL}}$ is monadic, and \Cref{CorBL} shows that the adjunction $F_{\mathrm{BL}}\dashv U_{\mathrm{BL}}$ is not an ideally exact context either. More generally, by \Cref{PropBLNotEquivalent}, no monadic functor $\BLAlg\to\BHoops$ has a left adjoint with cartesian unit. This answers the question raised in~\cite[Section~4]{MP}. A natural question, which we do not address here, is whether the ideally exact contexts associated with the unique morphisms from $\two$ to the terminal objects of $\GAlg$ and of $\BLAlg$ have a good theory of actions.

\section*{Acknowledgment}

The author is supported by the ``National Group for Algebraic and Geometric Structures and their Applications'' (GNSAGA -- INdAM), and he is a Postdoctoral Researcher of the Fonds de la Recherche Scientifique--FNRS.

%\subsection*{Competing interests}
%The author declares no competing interests.

\section*{Use of generative AI tools}
Generative AI tools (Anthropic's Claude Opus~5.5 and OpenAI's GPT-6~Sol) were used to improve the presentation of the content and to help check the correctness of the proofs. The author reviewed all suggestions made by these tools and takes full responsibility for the content of this article.

\vspace{1.5\baselineskip}

\end{document}